\documentclass[11pt]{amsart}
\usepackage{a4wide,tikz,hyperref,amsbsy,stmaryrd}

\usepackage{amsmath,amssymb,amsthm,latexsym,bm}
\usepackage{graphicx}
\usepackage{stmaryrd}		% \llbracket and \rrbracket
\usepackage{xcolor}
\usepackage{hyperref}
\usepackage{tikz}
\usepackage{amsbsy}

\usepackage{amssymb}
\usepackage{amsthm}
\usepackage{thmtools}
\usepackage{mathtools}
\usepackage{fullpage}

\usepackage[new]{old-arrows}

\usepackage[normalem]{ulem}

\usetikzlibrary{arrows,automata}
\usetikzlibrary{decorations.pathreplacing}
\usepackage[colorinlistoftodos]{todonotes}

\DeclareMathOperator{\Aut}{Aut}
\DeclareMathOperator{\supp}{supp}
\DeclareMathOperator{\Sym}{Sym}
\DeclareMathOperator{\Inn}{Inn}
\DeclareMathOperator{\Id}{Id}

\newtheorem{theorem}{Theorem}[section]
\newtheorem*{theorem*}{theorem}
\newtheorem{proposition}[theorem]{Proposition}
\newtheorem{lemma}[theorem]{Lemma}
\newtheorem{corollary}[theorem]{Corollary}

\theoremstyle{definition}   % <-- esto evita la itálica
\newtheorem{example}[theorem]{Example}
\newtheorem{rem}[theorem]{Remark}

\newtheorem{definition}[theorem]{Definition}

\def\Z{{\mathbb Z}}
\def\C{{\mathbb C}}
\def\N{{\mathbb N}}

\def\F{{\mathbb F}}

\def\vsp{{\vspace{0,2cm} }}

\title{Extensions of congruent monotileable (and amenable) groups
}

\author[P.~Cecchi Bernales]{Paulina Cecchi Bernales}
\address{Universidad de Chile, Departamento de Matem\'aticas, Santiago, Chile.}
\email{pcecchi@uchile.cl}
\author[E.~Molina]{Emir Molina}
\address{Universidad de Chile, Departamento de Matem\'aticas, Santiago, Chile.}
\email{emir.molina@ug.uchile.cl}
\author[C.~Rivas]{Crist\'obal Rivas}
\address{Universidad de Chile, Departamento de Matem\'aticas, Santiago, Chile.} 
\email{cristobalrivas@u.uchile.cl}
\date{\today}



\begin{document}

\begin{abstract}
    A monotile of a countable group $G$ is a finite subset $F \subseteq G$ that tiles $G$, namely $G=\bigsqcup_{c\in C} c F$, for some set of centers $C\subseteq G$. The group $G$ is called congruent monotileable (CM) if it admits an exhaustive sequence of monotiles $(F_n)_{n\in\N}$ such that each $F_{n+1}$ is  teselated by $F_n$. If moreover $G$ is amenable and the sequence $(F_n)_{n\in\N}$ is F{\o}lner, we say $G$ is congruent monotileable amenable (CMA). 

    In this work we observe that the classes of CM and CMA groups are stable under group  extensions satisfying a certain technical condition, which we call {\em stepwise monotileability}, and which is a weaker version of the notion of {\em $X$-monotileability} introduced by Dikranjan {\em et al.} Our condition also includes central extensions previously considered by Gao, Jackson, Seward for CM groups, and by Cecchi Bernales, Cortez for CMA groups. We show that this condition is satisfied by semidirect products of the form $G\wr H$, provided $G$ and $H$ are CM/CMA groups, and also by extensions of locally finite groups by CM/CMA groups. This furnishes new   examples of CM and CMA groups. In particular, we provide  the first examples of finitely generated CMA groups which are neither solvable nor residually finite.
\end{abstract}

\maketitle

\section{Introduction}
The study of monotileable groups was initiated by Chou in \cite{chou_80}. A (left) {\em monotile} $F$ of a countable group $G$ is a finite subset for which there exists a subset $\Delta\subseteq G$ such that $\{\delta F\}_{\delta\in \Delta}$ is a partition of $G$. The pair $(F,\Delta)$ is usually called a {\em tiling} of $G$ and $\Delta$ is called the set of {\em centers} of the tiling. 
The group is called {\em monotileable} (MT for short) if there exists an {\em exhaustive} sequence of monotiles, meaning that their union is the whole group $G$. 
Chou proved that every elementary amenable group is MT \cite{chou_80}. Later, Weiss raised the question of whether countable amenable groups admit F{\o}lner sets that are also monotiles \cite{weiss_2001}, and proved that residually finite amenable groups and solvable groups satisfy this property. The original relevance of monotiles comes from the fact that in amenable groups, a finite subset is a monotile if and only if it is a {\em Rokhlin set}, an object which plays a central role in the ergodic theory of amenable group actions \cite{ornstein_weiss_87}. More recently, Gao, Jackson and Seward introduced the concepts of ccc sequences of tilings and ccc groups in the context of studying marker structures on countable groups \cite{GJS}. Roughly speaking, a sequence of tilings $(F_n,\Delta_n)_{n\in\N}$ of a countable group $G$ is called {\em ccc (coherent, cofinal and centered)} if each monotile $F_n$ also tiles each $F_m$ for any $m\geq n$. This last property is the one captured by the word {\em coherent}. The group itself is called {\em ccc} if it admits a ccc sequence of tilings. 
In \cite{GJS} the authors prove that, inside the family of countable groups, the following classes are ccc: locally finite groups, residually finite groups, nilpotent groups, solvable groups with polycyclic commutator and free products of a collection of nontrivial groups. There are no known examples of countable groups which are not ccc. 
In another context, Cecchi Bernales and Cortez introduced the notion of {\em congruent monotileable amenable} groups (CMA for short) \cite{cecchi_cortez_2019}, which correspond to countable amenable groups admitting an exhaustive F\o lner sequence made of monotiles, which is also coherent in the sense of \cite{GJS}.

As noticed in \cite[Lemma 3.12]{salizzoni_et_al_23}, a congruent sequence necessarily consists of monotiles of the group (see also Proposition \ref{prop:equivalence} below). In \cite{cecchi_cortez_2019}, the authors prove that virtually nilpotent countable groups are CMA, and that for every CMA group $G$ and every Choquet simplex $K$ there exists a minimal $G$-subshit, which is free in a full measure set, whose set of invariant measures is affine homeomorphic to $K$. Moreover, Straszak proves in \cite{straszak_2020} that the class of minimal subshifts of CMA groups realizes any finite value of topological entropy, through what she calls {\em quasi-Toeplitz} subshifts. These last results illustrate the relevance of CMA groups in symbolic dynamics.

The class of CMA groups has also been considered by Delabie {\em et al.} in the context of measurable equivalence between amenable groups \cite{lemaitre_et_al_22}. More precisely, the authors use congruent F{\o}lner sequences of monotiles for constructing orbit equivalent couplings with prescribed integrability conditions between amenable groups (see \cite[Section 6]{lemaitre_et_al_22}). In the recent work  \cite{salizzoni_et_al_23}, Dikranjan {\em et al.} have extended the notion of CMA groups to monoids, under the name of {\em locally monotileable monoids}, in the context of studying algebraic entropy for endomorphisms of abelian groups. In that article, the authors prove that the notions of locally monotileable and CMA coincide for groups (\cite[Proposition 3.14]{salizzoni_et_al_23}) and provide a family of examples of CMA groups which are neither virtually nilpotent nor residually finite.

The notion of ccc groups previously discussed is equivalent to what we call in this article {\em congruent monotileable groups} (CM for short), following the vocabulary of \cite{cecchi_cortez_2019}. We prove the equivalence between both notions in Proposition \ref{prop:equivalence}. We remark that studying CM groups seems to be pertinent to the theory of hyperfinite orbit equivalence relations. Indeed, as pointed out in \cite{GJS}, CM sequences of monotiles constitute a strong type of marker structures, a tool that has been used several times in studying hyperfinite equivalence relations (we refer to \cite[Chapter 4]{GJS} for details). On the other hand, the question of whether the orbit equivalence relation generated by a Borel action of any countable amenable group is hyperfinite, first proposed by Weiss in \cite{weiss_measurable_dynamics:1984}, remains open. This suggest a natural connection between CMA groups and hyperfinite orbit equivalence relations. More recently, Bitar, Cabezas and Guillon \cite{BCG26} have generalized the notions of substitutions and $S$-adic systems from the classical $\Z$ systems to actions of general countable groups by using congruent sequences of (possibly multiple) tiles. They call such sequences {\em hierarchical tiling structures}, and use them to somehow emulate the linear geometry of $\Z$ in order to define substitutions and directive sequences of morphisms. The particular cases of CM and CMA groups are respectively called {\em uniformly hierarchical} and {\em F{\o}lner hierarchical} groups in \cite{BCG26}.

It is an open question whether every amenable group is CMA. Note in addition that, {\em a priori}, a countable amenable CM group is not necessarily CMA, since being CMA requires the congruent sequence of monotiles to be F{\o}lner. In fact, it is an open question whether any CM group which is also amenable is CMA. 

\vsp
In this article, we study extension properties of CM and CMA groups. Consider a short exact sequence
\begin{equation} \label{eq:shortexact}
 1\longrightarrow K\longrightarrow  G\longrightarrow H\longrightarrow 1.   
\end{equation}  where $K$, $G$ and $H$ are countable groups. 
We always assume that $K$ is a subgroup of $G$.
Our goal is to find sufficient conditions on the CM/CMA groups $K$ and $H$ so that the CM/CMA property is inherited by the group $G$. This problem has been addressed in \cite{cecchi_cortez_2019} for CMA groups, where the authors prove that a sufficient condition is that $K$ be a central subgroup of $G$. An analogous condition is obtained in \cite{GJS} in the context of CM groups. 
In \cite{salizzoni_et_al_23}, the authors introduce the following  notion: for $X\subseteq \Aut(G)$, a sequence $(F_n)_{n\in \N}$ of finite subsets of $G$ is said to be {\em $X$-monotileable} if for all $n\in \N$ and $\phi\in X$, the set $F_{n+1}$ is tiled by $\phi(F_{n})$, {\em i.e.}, there exists a subset $J_{n,\phi}\subseteq G$ such that $F_{n+1}=\bigsqcup_{j\in J_{n,\phi}}j\phi(F_n)$. The group $G$ is called $X$-monotileable if it admits an $X$-monotileable F\o lner sequence.
In  \cite[Theorem 1.15]{salizzoni_et_al_23}, the authors prove that  $G$ is CMA whenever $K, H$ are CMA and $K$ is $\Inn(G)$-monotileable, where $\Inn(G)$ denotes the set of all inner automorphisms of $G$, or, equivalently, the set of all conjugation maps. 

Our approach to the previous problem consists in using a weaker notion of $X$-monotileability. Let $(X_n)_{n\in \N}$ be a sequence of subsets of $\Aut(G)$. A sequence $(F_n)_{n\in \N}$ of finite subsets of $G$ is called {\em $(X_n)$-stepwise congruent} if for each $n\in\N$ and each $\phi\in X_n$, the set $F_{n+1}$ is tiled by $\phi(F_{n})$. 
If $K, G$ and $H$ are countable groups like in the short exact sequence \eqref{eq:shortexact}, where $K$ and $H$ are CM/CMA groups, we say that $K$ is {\em $G$-stepwise monotileable} if there exists a congruent sequence $(F_n)_{n\in\N}$ of finite subsets of $K$, a sequence of liftings $(\widehat{T}_n)$ of a congruent sequence $(T_n)_{n\in\N}$ of $H\cong G/K$, and a sequence $(\widehat{S}_n)_{n\in\N}$, such that 
\begin{itemize}
    \item $(\widehat{S}_n)_{n\in\N}$ is a sequence of {\em local centers} of $(\widehat{T}_n)_{n\in\N}$, that is, 
    $\widehat{T}_{n+1}=\bigsqcup_{\hat{s}\in \widehat{S}_n} \hat{s}\widehat{T}_n$ for all $n$, and
    \item $(F_{n})_{n\in\N}$ is $(\widehat{S}_n)$-stepwise congruent.
\end{itemize}

We note in Proposition \ref{prop:innimpliesstepwise} that $\Inn(G)$-monotileability implies $G$-stepwise monotileability. Using this notion we show the main abstract result of this paper, which is a condition that ensure congruent gluing of the tiles in a short exact sequence.
\begin{theorem}\label{thm:main}
    Let $H,G,K$ be countable groups. Suppose there exists a short exact sequence  
    $$1\longrightarrow K\longrightarrow  G\longrightarrow H\longrightarrow 1,$$
    and $K$, $H$ are both CM. If $K$ is $G$-stepwise monotileable, then $G$ is a CM group.
\end{theorem}
In the setting of CMA groups, we further prove in Theorem \ref{thm:folner} that under stepwise monotileability, if the sequences $(F_m)_{m\in\N}$ and $(T_n)_{n\in\N}$ are F{\o}lner, and an increasing sequence $(m_n)_{n\in\N}$ is carefully chosen so that $(F_{m_n}\widehat{T}_n)_{n\in\N}$ is F{\o}lner in $G$ (see Lemma \ref{lem:subsequencefolner}), then CMA of $K$ and $H$ is inherited by $G$.\\

Applying Theorems \ref{thm:main} and \ref{thm:folner} we identify new stability properties for the class of CM/CMA groups. 

\begin{theorem} \label{thm:locallyfinitebycma}
Let $H$, $G$, $K$ be countable groups. Suppose there exists a short exact sequence 
  $$1\longrightarrow K\longrightarrow  G\longrightarrow H\longrightarrow 1,$$
$K$ is locally finite and $H$ is CM/CMA. Then, $G$ is CM/CMA. 
\end{theorem}

\begin{theorem} \label{thm:wr}
Let $G,H$ be both CM/CMA groups and let $X$ be a countable set on which $H$ is acting by  permutations. Then permutational wreath product $G\wr_X H$ is also a CM/CMA group.
\end{theorem}

The particular case of Theorem \ref{thm:wr} for $G=(\Z/m\Z)\wr \Z$, $m\geq2$, was proved in \cite{lemaitre_et_al_22}. However, Theorem \ref{thm:wr} implies more generally that the class of CMA groups is stable under wreath product, thus getting new  examples of CMA groups which are not residually finite (for instance $H\wr \Z$ for an $H$ non-abelian, see Remark \ref{rem:Grunberg}). Furthermore, from Theorem \ref{thm:locallyfinitebycma} we obtain examples of finitely generated CMA groups which are neither residually finite, nor virtually solvable. For instance, $G=\Sym_0(\Z)\rtimes\Z$, where $\Sym_0(\Z)$ denotes the group of finitely supported permutations of $\Z$. To the best of our knowledge, these constitute the first examples of their kind.\\

The text is organized as follows. In Section \ref{sec:tilings}, we give all necessary preliminaries to this work, including the precise definitions of congruent monotileable and congruent monotileable amenable groups. We review and compare the notion of congruent monotileable groups with the previous one of ccc groups introduced in \cite{GJS}. 
In Section \ref{sec:mainresults}, we present the main results of the article. In particular, we give the proof of Theorems \ref{thm:main} and  \ref{thm:folner} and compare our notion of stepwise monotileability with the previous one of $X$-monotileability \cite{salizzoni_et_al_23}.
All the examples of groups studied in this paper are presented in Section \ref{sec:examples}. In this section, we prove Theorems \ref{thm:locallyfinitebycma} and \ref{thm:wr}, and specify some families of groups for which these theorems apply. Finally, as a last consequence of Theorems \ref{thm:main} and \ref{thm:folner}, we treat the example of the {\em Lampdesigner group} $(G\wr_H \Sym_0(H))\rtimes H$, and prove that this group is CM/CMA provided that $G$ and $H$ are CM/CMA as well. \\

Throughout the text, groups are countable unless otherwise stated. Given a group $G$, the identity element of $G$ is denoted $1_G$.

\bigskip

\paragraph{{\bf Acknowledgements}}  
We would like to thank Mar\'ia Isabel Cortez for suggesting that stepwise monotileability could be used to prove Theorem \ref{thm:locallyfinitebycma}, as well as Michele Triestino and Mart\'in Gilabert for bringing Examples \ref{ex:houghton} and \ref{ex:lampdesigner} to our attention. The authors acknowledge the support of Math Amsud project AMSUD250017. P. Cecchi Bernales and E. Molina were supported by ANID/Fondecyt Iniciaci\'on 11240886. C. Rivas was supported by ANID/Fondecyt 1241135.
\subsection*{AI disclosure}
No AI engine were used for the preparation of this manuscript: every mathematical
claim, grammar style and typo is of human nature.

\section{Tilings of Countable Groups}\label{sec:tilings}

In this section we introduce the main combinatorial objects used in this work, namely congruent sequences of monotiles, which are the base for the definitions of congruent monotileable (CM) and congruent monotileable amenable (CMA) groups (see Definitions \ref{def:congruent} and \ref{def:cma}). We first compare the notion of CM groups with that of {\em ccc groups} introduced in \cite{GJS} and  prove that they are equivalent in Proposition \ref{prop:equivalence}. Then in Section \ref{subsec:CMA} we treat the particular case of CMA groups.

Let $G$ be a group, let $F, V \subseteq G$. We say that $F$ {\em tiles} $V$ if there exists a subset $C\subset G$ such that $V=\bigsqcup_{c\in C} cF$. In the case $V=G$, we call $F$ a {\em monotile} of $G$ and the pair $(F,C)$ a {\em tiling} of $G$. If a finite subset $F\subseteq G$ verifies that $G=\bigsqcup_{c\in C} Fc$ for some $C\subseteq G$, $F$ is called a {\em right monotile} of $G$.\\
A sequence of tilings $(F_n,\Delta_n)_{n\in \N}$ is called \textit{coherent} if whenever $n\leq m$ and $\delta_n\in \Delta_n$, $\delta_m\in \Delta_m$ satisfy $\delta_nF_n\cap \delta_mF_m\neq \emptyset$, then $\delta_nF_n\subseteq \delta_mF_m$. 
It is called \textit{cofinal} if for every finite subset $A\subset G$ there exist $n\geq 1$ and $\delta\in \Delta_n$ such that $A\subseteq \delta F_n$, and \textit{centered} if $1_G\in \Delta_n$ for every $n\geq 1$. 
We say that the sequence $(F_n, \Delta_n)_{n\in\N}$ is a {\em ccc tiling} if it is coherent, cofinal and centered.  A countable group $G$ is called {\em ccc} if it admits a ccc tiling. The definitions of ccc tilings and ccc groups were introduced in \cite{GJS}.  \\

The following combinatorial objects were introduced in \cite{cecchi_cortez_2019}. Although they were originally defined for countable amenable groups, here we first consider the more general version that does not depend on the amenability of the group. 

\begin{definition}\label{def:congruent}
   Let $G$ be a group. A sequence $(F_n)_{n\in\N}$ of finite subsets of $G$ is called congruent if $1_G\in F_0$ and for every $n\geq 0$ there exists a set $J_n\subset G$ such that $1_G\in J_n$ and $F_{n+1}=\bigsqcup_{j\in J_n}jF_n$. We call  $(J_n)_{n\in \N}$ a sequence of local centers associated to $(F_n)_{n\in \N}$. The group $G$ is called congruent monotileable (CM) if it admits a congruent and exhaustive sequence of monotiles.
\end{definition}

By definition, a congruent exhaustive sequence $(F_n)_{n\in\N}$ satisfies the following conditions,
\begin{itemize}
        \item[(i)] $1_G\in F_0$
        \item[(ii)] For every $n\geq 0$, there exists a finite subset $J_n\subseteq G$ such that $1_G\in J_n$ and $F_{n+1}= \bigsqcup_{j\in J_n} jF_n$.
        \item[(iii)] $G=\bigcup_{n\geq 0} F_n$.
\end{itemize}

We prove that these three conditions force each $F_n$ to be a monotile. This was also noticed in \cite[Lemma 3.12]{salizzoni_et_al_23}. In addition, we prove that conditions (i)-(iii) characterize both the ccc tilings and congruent exhaustive sequences, showing that both combinatorial structures are equivalent.

\begin{proposition} \label{prop:equivalence}
Let $G$ be a countable group. The following are equivalent,
\item[1)] $G$ is ccc,
\item[2)] $G$ is CM, and
\item[3)] $G$  admits a sequence of finite subsets that satisfies conditions (i)-(iii).
\end{proposition}
\begin{proof}
$2) \Longrightarrow 3)$ is clear. We prove $1)\Longrightarrow 2)$ and $3) \Longrightarrow 1)$.\\

$1)\Longrightarrow 2)$ Suppose $G$ admits a ccc tiling $(F_n, \Delta_n)_{n\in\N}$. 
Note that we may always assume that $1_G\in F_0$ by replacing all $\Delta_n$ and $F_n$ with $\Delta_ng$ and $g^{-1}F_n$, respectively, for some $g\in F_0$.
Moreover, since $(F_n, \Delta_n)_{n\in\N}$ is coherent and $G=\bigsqcup_{\delta\in \Delta_n}\delta F_n$ for all $n\geq 0$, we have that for all $n\geq 0$ there exists a finite subset $J_n\subset \Delta_n$ such that $F_{n+1}=\bigsqcup_{j\in J_n}j F_n$. This means that $(J_n)_{n\in \N}$ is a sequence of local centers associated to $(F_n)_{n\in\N}$.
On the other hand, since $(F_n, \Delta_n)_{n\in\N}$ is cofinal, for every finite subset $A\subset G$ there exist $n\geq 0$ and $\delta \in \Delta_n$ such that $A\cup\{1_G\}\subseteq \delta F_n$. Since $1_G\in \delta F_n\cap F_n$, we deduce that $\delta=1_G$. Hence, $A\subset A\cup\{1_G\}\subseteq F_n$. This implies that $G=\bigcup_{n\geq1}F_n$. We conclude that $(F_n)_{n\in \N}$ is congruent and exhaustive. Therefore, $G$ is a CM group.\\

$3) \Longrightarrow 1)$ Let $(F_n)_{n\in\N}$ be a sequence of finite subsets of $G$ satisfying properties (i)-(iii), let $(J_n)_{n\in\N}$ be the associated sequence of local centers.
For every $n\geq 0$, define $\Delta_n:=\bigcup_{m\geq n}J_m\cdots J_n$. From conditions (ii) and (iii), it follows that $G=\bigcup_{\delta \in \Delta_n}\delta F_n$ for every $n\geq0$.
Observe that for each $n\geq 0$, for every $m>n$, the collection
$$\{j_{m-1}\cdots j_nF_n\ :\ j_i\in J_i,\ \forall n\leq i\leq m-1 \}$$
is a partition of $F_m$. Indeed, it suffices to apply condition (ii) inductively. From this fact we verify that $\delta F_n\cap\delta'F_n=\emptyset$ whenever $\delta, \delta'\in \Delta_n$ and $\delta\neq \delta'$. 

Consider two distinct elements $\delta=j_m\cdots j_n$, and $\delta=j_{\tilde{m}}'\cdots j_n'\in \Delta_n$.
If $m=\tilde{m}$, we may suppose without loss of generality that $j_m\neq j'_{m}$. Hence, the fact that $j_mF_m\cap j'_{m}F_m=\emptyset$ implies that $\delta F_n\cap \delta' F_n=\emptyset$. 
On the other hand, if $m\neq \tilde{m}$ with $m<\tilde{m}$, then $F_{m+1}\cap j'_{\tilde{m}}\cdots j'_{m+1}F_{m+1}=\emptyset$. Hence, $\delta F_n\cap \delta'F_n=\emptyset$.
Finally, we have $G=\bigsqcup_{\delta \in \Delta_n}\delta F_n$.
Therefore, every $(F_n, \Delta_n)$ is a tiling of $G$.

We now prove that $(\Delta_n, F_n)_{n\in\N}$ is coherent, cofinal and centered. Since $1_G\in J_n$ for every $n\geq 0$, every pair $(F_n, \Delta_n)$ is a centered tiling. Moreover, by the definition of $\Delta_n$, condition (iii) immediately implies that $(F_n, \Delta_n)_{n\in\N}$ is cofinal. In addition, note that $(F_n, \Delta_n)_{n\in\N}$ is coherent. Indeed, for all $n\leq m$ and $\delta_m\in \Delta_m$ 
$$\delta_mF_m=\bigsqcup_{j_{n}\in J_{n},\dots, j_{m-1}\in J_{m-1}} \delta_mj_{m-1}\ldots j_nF_n.$$

We conclude that $(F_n, \Delta_n)_{n\in\N}$ is a ccc tiling and, therefore, $G$ is a ccc group.
\end{proof}

\begin{rem}
Proving that every element in a sequence of finite subsets is a monotile might become tedious work, especially when the choice of such a sets is not ``natural''. Proposition \ref{prop:equivalence} shows that this work might be avoided by verifying conditions (ii) and (iii).
In each of the proofs below, we will limit ourselves to verify conditions (i), (ii), (iii) every time we want to prove that a sequence of finite subsets is a congruent and exhaustive sequence of monotiles.
\end{rem}

\subsection{Congruent Monotileable Amenable Groups} \label{subsec:CMA}

Recall that a sequence  $(F_n)_{n\in \N}$ of finite subset of a group $G$ is called a right F\o lner sequence if for every $g\in G$, 
\[ \lim_{n\to\infty}\frac{|F_n\setminus F_ng|}{|F_n|}=0.\]
Equivalently, $(F_n)_{n\in \N}$ is a right F\o lner sequence if for every finite subset $C\subseteq G$ and every $\varepsilon>0$ there exists $N\in \N$ such that for all $n\geq N$ the set $F_{n}$ is right $(C,\varepsilon)$-invariant, that is
\[|\{g\in F_n\ |\ gC\subseteq F_n\}| \geq (1-\varepsilon)\cdot |F_n|.\]
Left F{\o}lner sequences and left $(C,\varepsilon)$-invariance are defined analogously. In this article we only work with right F{\o}lner sequences, so unless otherwise stated, {\em F{\o}lner sequence} means {\em right F{\o}lner sequence}.\\
A group that admits a F\o lner sequence is called {\em amenable}.

\begin{definition}\label{def:cma}
    Let $G$ an amenable group. We say that $G$ is congruent monotileable amenable (CMA) if it admits a congruent and exhaustive F\o lner sequence.
\end{definition}

\begin{rem}\label{rem:increasing}
Every countable amenable group $G$ admits an exhaustive and increasing F\o lner sequence $(F_n)_{n\in \N}$ such that every set $F_n$ contains the identity of the group. The property of $(F_n)_{n\in\N}$ being congruent is more delicate. There are no known examples of amenable groups which are not CMA, but the question of whether every amenable group is CMA remains open.
\end{rem}

It is known that countable amenable groups which are also residually finite \cite{cortez_petite_2014}, or virtually nilpotent \cite{cecchi_cortez_2019}, or locally finite \cite{salizzoni_et_al_23} are all CMA.
In \cite{GJS}, Gao, Jackson, and Seward proved that countable locally finite groups are all CM \cite[Theorem 4.5.3]{GJS}. A short proof of this fact is provided in \cite{salizzoni_et_al_23}. In order to show how our characterization with conditions (i), (ii), (iii) simplifies the details of the proof given in \cite{GJS}, we provide yet another proof in the following.

\begin{proposition}\label{prop:locallyfinite}
 Every countable locally finite group is CMA.   
\end{proposition} 
\begin{proof}
    Let $G$ be a countable locally finite group and $(A_n)_{n\in \N}$ a sequence of finite subgroups of $G$ such that $G=\bigcup_{n\in \N}A_n$ and $A_n\lneq A_{n+1}$ for all $n\in \N$.     
    Since every $A_n$ is finite, it suffices to consider $J_n\subset A_{n+1}$ a set of coset representatives of $A_{n+1}/A_n$. Then, $(J_n)_{n\in \N}$ is a sequence of local centers associated to $(A_n)_{n\in \N}$ and therefore $(A_n)_{n\in \N}$ is congruent.
    Moreover, since $G=\bigcup_{n\in N} A_n$, the sequence $(A_n)_{n\in \N}$ is also F\o lner. 
    Indeed, if $g$ is an arbitrary element of $G$, then there exists $n_0\in \N $ such that $g\in A_{n_0}$. Consequently, $|A_n\setminus A_ng|=|A_n\setminus A_n|=0,\ \forall n\geq n_0.$
\end{proof}

Let $G,K,H$ be countable groups satisfying a short exact sequence as in Equation \eqref{eq:shortexact}, with  $p\colon G\to H\cong G/K$ the canonical projection. Given any set $T\subseteq H$, we call $\widehat{T}\subseteq G$ a {\em lifting} of $T$ if $p$ induces a bijection from $ \widehat{T}$ onto $T$.

The following result due to Weiss \cite{weiss_2001} constitute the key idea to prove that admitting F{\o}lner sequences of monotiles is stable under group extensions. We will use it later to prove Lemma \ref{lem:subsequencefolner} which in turn is needed for the proof of Theorems \ref{thm:locallyfinitebycma}, \ref{thm:wr} and \ref{thm:folner}. Remark that the original version in \cite{weiss_2001} is stated for left invariance and right monotiles, while here we use right invariance and left monotiles, however, the proof is completely analogous.

\begin{theorem}\cite[Theorem 2]{weiss_2001}\label{thm:weisslemma}
    Let $G$ be a countable group and let $K$ be a normal subgroup of $G$. Let $p \colon G\to G/K\cong H$ be the canonical homomorphism from $G$ onto $G/K$, $C$ a finite subset of $G$ and $\varepsilon>0$. If $T$ is a monotile in $H$ which is $(p(C),\varepsilon/2)$-invariant and $\widehat{T}$ is a  lifting of $T$, then there exists $J\subseteq K$ such that for any monotile $U\subseteq K$ which is $(J,\varepsilon/2)$-invariant, the set $U\widehat{T}$ is a $(C,\varepsilon)$-invariant monotile of $G$.
\end{theorem}
 The following lemma is very similar to \cite[Lemma 7]{cecchi_cortez_2019}, adapted to our situation to get item (2) of the statement. Its proof strongly relies on Theorem \ref{thm:weisslemma}.

\begin{lemma} \label{lem:subsequencefolner}
    Let $G,K,H$ be countable groups that satisfy a short exact sequence
        $$1\longrightarrow K\longrightarrow  G\longrightarrow H\longrightarrow 1.$$
    Suppose that $K$ and $H$ have 
    exhaustive F\o lner sequences of monotiles $(F_n)_{n\in \N}$ and $(T_n)_{n\in\N}$, respectively. Let $(R_n)_{n\in\N}$ be any sequence of finite subsets of $G$. Then, for any sequence of liftings $(\widehat{T}_n)_{n\in \N}$,
    there exists an increasing sequence $(m_n)_{n\in \N}$ such that 
    \begin{enumerate}
        \item $\widetilde{F}_n:=F_{m_n}\widehat{T}_n$ defines an exhaustive F\o lner sequence of monotiles of $G$, and
        \item  for all $n\in\N$, $\bigcup_{r\in R_n} rF_{m_n}r^{-1}\subseteq F_{m_{n+1}}$.
    \end{enumerate}
     \end{lemma}

\begin{proof}
Let $p\colon G\to H$ be the canonical projection. Let $(C_n)_{n\in\N}$ be an increasing sequence of finite subsets of $G$ such that $G=\bigcup_{n\in\N}C_n$, let $(\varepsilon_n)_{n\in\N}$ be a decreasing to zero sequence of positive real numbers. Since $(T_n)_{n\in\N}$ is a F{\o}lner sequence we have that, up to taking a subsequence, each $T_n$ is right $(p(C_n),\varepsilon_n/2)$-invariant. We can also assume that $T_0=\{1_H\}$.\\ We inductively define $(m_n)_{n\in\N}$. Let  $m_0=0$. 
Let $n>0$. Suppose we have defined 
$m_{n-1}$. 
Since both $R_{n-1}$ and $F_{m_{n-1}}$ are finite sets, and the sequence $(F_m)_{m\in\N}$ is exhaustive, there exists $M_1\in \N$ such that for all $m\geq M_1$, $\bigcup_{r\in R_n} rF_{m_{n-1}}r^{-1}\subseteq F_m$.\\
On the other hand, since $T_n$ is a right $(p(C_n),\varepsilon_n/2)$-invariant monotile of $H$, by Theorem \ref{thm:weisslemma}, there exists a finite subset $J\subseteq K$ such that if $F\subseteq K$ is a $(J,\varepsilon_n/2)$-invariant monotile of $K$, then $F\widehat{T}_n$ is a $(C_n,\varepsilon_n)$-invariant right monotile of $G$. Since $(F_m)_{m\in\N}$ is a F{\o}lner sequence, there exists $M_2\in\N$ such that for all $m\geq M_2$, $F_m$ is right $(J,\varepsilon_n/2)$-invariant. Choose $m_n>\max\{M_1,M_2\}$. We get that $F_{m_n}\widehat{T}_n$ is a right $(C_n,\varepsilon_n)$-invariant monotile of $G$. This shows that $(F_{m_n}\widehat{T}_n)_{n\in\N}$ is a F{\o}lner sequence made of monotiles of $G$. Moreover, the fact that $(F_{m_n}\widehat{T}_n)_{n\in\N}$ is exhaustive follows directly from the fact that $(F_n)_{n\in\N}$ and $(T_n)_{n\in\N}$ are. This proves (1). By construction, $\bigcup_{r\in R_n} rF_{m_n}r^{-1}\subseteq F_{m_{n+1}}$ for all $n$, which proves (2) and concludes the proof of the lemma. 
\end{proof}

\section{Main results}\label{sec:mainresults}

Let $G,K,H$ be countable groups satisfying a short exact sequence as in Equation \eqref{eq:shortexact}.
If $(T_n)_{n\in \N}$  is a congruent and exhaustive sequence of $H$, with associated sequence of local centers $(S_n)_{n\in \N}$, then for any lifting of $(S_n)$ containing the identity we can inductively define a sequence of liftings $(\widehat{T}_n)_{n\in \N}$ such that $(\widehat{T}_n)$ is congruent with local centers $(\widehat{S}_n)$. For this we set $\widehat{T}_0$ any lifting of $T_0$, and $(\widehat{S}_n)$  liftings of $(S_n)$ containing the identity. Then $\widehat{T}_{n+1}:=\bigsqcup_{\hat{s}\in \widehat{S}_n}\hat{s}\widehat{T}_n$ defines a lifting of $T_{n+1}$  (see \cite[Lemma 7]{cecchi_cortez_2019} for details). If we now consider a congruent exhaustive sequence $(F_n)_{n\in \N}$ in $K$, we might expect that the natural candidate $(\widetilde{F}_n:=F_n\widehat{T}_n)_{n\in \N}$ preserves the same structure (i.e. is congruent). However this is not always true (for instance this fails already for some extension $\Z^2\rtimes \Z$ such as Sol or the Heisenberg group).
What is known is that the conclusion holds if we assume that $K$ is contained in the center of $G$, a  result proved  independently by Gao, Jackson, Seward \cite[Lemma 4.4.2]{GJS} and Cecchi Bernales, Cortez \cite[Lemma 7] {cecchi_cortez_2019}.
Another sufficient condition to conclude that $(\widetilde{F}_{n})_{n\in \N}$ is congruent and exhaustive is given by the notion of $X$-monotileability introduced in \cite{salizzoni_et_al_23} and recalled in the Introduction. This notion is useful, for instance, when there is a rigid subgroup structure, like in the following example.

\begin{example}[The Pr\"ufer group]\label{ex:prufer}
    Let $p\in \N$ be a prime number. The Pr\"ufer group $\Z(p^\infty)$ is defined as the following subgroup of the unit circle,
$$\Z(p^\infty)=\{\exp(i2\pi m/p^n):0\leq m<p^n,n\in\N\},$$
that is, the set of all $p^n$-th roots of unity, when $n$ runs over $\N$. Equivalently, $\Z(p^\infty)$ can be represented as the direct limit
$$\Z(p^\infty)=\varinjlim (\Z/p^n\Z,i_n),$$
where $i_n:\Z/p^n\Z\to \Z/p^{n+1}\Z$ is the multiplication by $p$. We refer to \cite{Fu} for details on the properties of $\Z(p^\infty)$. The Pr\"ufer group is countable, locally finite and divisible, thus not residually finite. Moreover, for a given $n\in\N$, there is  only one subgroup of order $p^n$ in $\Z(p^\infty)$, namely  $\Z/p^n\Z$. Note that the sequence $(A_n)_{n\in\N}$ defined by $A_n=\Z/p^n\Z$ is an increasing exhaustive sequence of finite subgroups. As explained in Proposition \ref{prop:locallyfinite}, $(A_n)_{n\in\N}$ is also congruent and F\o lner. 

Now, if $G$ is a group containing $\Z(p^\infty)$ as a normal subgroup and $\phi\in \Inn(G)$, then $\phi(A_n)$ has order $p^n$ and is a subgroup of $\Z(p^\infty)$. Thus by uniqueness it necessarily equals $A_n$. This proves that $\Z(p^\infty)$ is $\Inn(G)$-monotileable. It follows from \cite[Theorem 1.15]{salizzoni_et_al_23} that if $G/\Z(p^\infty)$ is CMA, so it is $G$. 

\end{example}

Our main observation is that one can still deduce CM/CMA for group extensions even in situations when there is a much less rigid subgroup structure. For this we introduce  the notions of {\em stepwise congruence} and {\em stepwise monotileability}, which constitute the key condition under which we can prove that extensions of CM/CMA groups are also CM/CMA.

\begin{definition}\label{def:stepwisecongruent}
    Let $G$ be a group, let $(X_n)_{n\in \N}$ be a sequence of subsets of $\Aut(G)$. A sequence $(F_n)_{n\in \N}$ of finite subsets of $G$ is said to be $(X_n)$-stepwise congruent if for each $n\in\N$ and each $\phi\in X_n$, the set $F_{n+1}$ is tiled by $\phi(F_{n})$.
\end{definition}

\begin{definition}\label{def:stepwisemonotileable}
 Let $K, G$ and $H$ be countable groups like in the short exact sequence \eqref{eq:shortexact}, where $K$ and $H$ are CM/CMA groups. We say that $K$ is $G$-stepwise monotileable if there exists a congruent and exhaustive sequence $(F_n)_{n\in\N}$ of $K$, a congruent and exhasutive sequence $(T_n)_{n\in\N}$ of $H\cong G/K$ together with a sequence of local centers $(S_n)_{n\in\N}$, and liftings $(\widehat{T}_n)_{n\in\N}$, $(\widehat{S}_n)_{n\in\N}$, such that $(\widehat{S}_n)_{n\in\N}$ is a sequence of local centers of $(\widehat{T}_n)_{n\in\N}$, and $(F_{n})_{n\in\N}$ is $(\widehat{S}_n)$-stepwise congruent.   
\end{definition}

Although the notions of $\Inn(G)$-monotileability and stepwise monotileability might look similar, they are not the same. In the following proposition we note that $\Inn(G)$-monotileability implies stepwise monotileability. On the other hand, we prove in Proposition \ref{prop:wreath_notInnmonotileable} and Proposition \ref{prop:sym_notInnmonotileable} that stepwise monotileability does not necessarily imply $\Inn(G)$-monotileability, which shows that $\Inn(G)$-monotileability is strictly stronger than stepwise monotileability.

\begin{proposition}\label{prop:innimpliesstepwise}
    Let $K,G,H$ be countable groups satisfying a short exact sequence as in Equation \eqref{eq:shortexact}. Suppose that $K$ and $H$ are CM/CMA. If $K$ is $\Inn(G)$-monotileable, then $K$ is $G$-stepwise monotileable.
\end{proposition}
\begin{proof}
    Let $(F_n)_{n\in\N}$ be a congruent exhaustive sequence of finite subsets of $K$ which is also $\Inn(G)$-monotileable, and let $(T_n)_{n\in \N}$ be any congruent exhaustive sequence of finite subsets of $H$, together with centers sequence $(S_n)_{n\in\N}$. Consider the liftings $(\widehat{T}_n)_{n\in \N}$ and $(\widehat{S}_n)_{n\in \N}$ as before. For every $\hat{s}\in \widehat{S}_n$ and every $n\in \N$, there exists an inner automorphism $\varphi_{\hat{s}}$ defined by $\varphi(g):=\hat{s}g\hat{s}^{-1}$.
Since $(F_n)_{n\in \N}$ is $\Inn(G)$-monotileable, $\varphi_{\hat{s}}(F_n)=\hat{s}F_{n}\hat{s}^{-1}$ tiles $F_{n+1}$ for every $\hat{s}\in \widehat{S}_n$, for every $n\in \N$. In particular, $(F_n)_{n\in \N}$ is also a $(\widehat{S}_n)$-stepwise congruent sequence. This proves that $K$ is $G$-stepwise monotileable. 
\end{proof}

We now prove the main abstract result of the article.
\vspace{0,1cm}

\noindent\textit{Proof of Theorem \ref{thm:main}.}
Since $K$ is $G$-stepwise monotileable, there exists a congruent exhaustive sequence $(F_n)_{n\in\N}$ of finite subsets of $K$, a congruent exhaustive sequence $(T_n)_{n\in\N}$ of $H$ together with a sequence $(S_n)_{n\in\N}$ of local centers and liftings $(\widehat{T}_n)_{n\in\N}$, $(\widehat{S}_n)_{n\in\N}$, such that $(\widehat{S}_n)_{n\in\N}$ is a sequence of local centers of $(\widehat{T}_n)_{n\in\N}$, and $(F_{n})_{n\in\N}$ is $(\widehat{S}_n)$-stepwise congruent. We will prove that $(\widetilde{F}_n:=F_n\widehat{T}_n)_{n\in\N}$ is a congruent and exhaustive sequence of finite subsets of $G$. According to Proposition \ref{prop:equivalence}, it suffices to show that this sequence satisfies conditions (i)-(iii) from Section \ref{sec:tilings}.\\
First,  note that $1_G\in \widetilde{F}_0$ since $1_G\in F_0\cap \widehat{T}_0$. Thus, condition (i) follows immediately.\\
Consider now an arbitrary element $g\in G$ and a nonnegative integer $n_0$ such that $p(g)\in T_{n_0}$, where $p\colon G\to H$ denotes the canonical projection.
Since $T_{n_0}=S_{n_0-1}T_{n_0-1}=\dots=S_{n_0-1}\dots S_0T_0$, there exist an element $t\in T_0$ and elements $s_i\in S_i$ for every $i\in\{ 0,\dots, n_0-1\}$ such  that $p(g)=s_{n_0-1} \dots s_0t$. By the definition of $\widehat{T}_0$ and the $\widehat{S}_n$'s, this means that $g\in K\hat{s}_{{n_0}-1}\dots \hat{s}_{{1}}\hat{t}$.
Therefore, $g= k\hat{s}_{{n_0}-1}\dots \hat{s}_{{0}}\hat{t}$ for some $k\in K$.

On one hand, $\hat{s}_{n_0-1}\dots \hat{s}_{{0}}\hat{t}\in \widehat{S}_{n_0-1}\dots \widehat{S}_0\widehat{T}_0=\widehat{T}_{n_0}$.
On the other hand, since $(F_m)_{m\in\N}$ is exhaustive, $k\in F_{m_0}$ for some $m_0\in \N$. Thus, taking $N=\text{max}\{n_0,m_0\}$, we have
$$g=k\hat{s}_{{n_0}}\dots \hat{s}_{{1}}\hat{t}\in F_{m_0}\widehat{T}_{n_0}\subseteq F_N\widehat{T}_N=\widetilde{F}_N.$$
This proves that $(\widetilde{F}_n)_{n\in \N}$ satisfies (iii), {\em i.e,}, that $(\widetilde{F}_n)_{n\in \N}$ is exhaustive. \\
Finally, let us show that every set $\widetilde{F}_{n+1}$ is the union of pairwise disjoint translations of $\widetilde{F}_n$, which will prove condition (ii).

Let $n\in \N$. Consider an element $s\in S_n$ and the finite subset $\Delta_{n}^{(s)}\subseteq G$ such that $F_{n+1} = \bigsqcup_{\delta\in \Delta_{n}^{(s)}}\delta (\hat{s}F_n\hat{s}^{-1})$, which exists by the hypothesis of stepwise monotileability. After a simple rewriting, we get that
$$\Delta_n^{(s)}\hat{s}\cdot \widetilde{F}_n=\Delta_n^{(s)}\hat{s}\cdot F_n\widehat{T}_n=\Delta_{n}^{(s)} (\hat{s}F_n\hat{s}^{-1})\cdot \hat{s}\widehat{T}_n=F_{n+1}\cdot \hat{s}\widehat{T}_n.$$

Therefore, the set $\widetilde{J}_n=\bigcup_{s\in S_n}\Delta_{n}^{(s)}\hat{s}$ is finite, contains the identity (by construction of the sets $\Delta_n^{(s)}$'s and $\widehat{S}_n$'s) and, moreover, it satisfies 
\begin{align*}
    \widetilde{F}_{n+1}= F_{n+1}\widehat{T}_{n+1}&= \bigcup_{s\in S_n} F_{n+1} \cdot\hat{s}\widehat{T}_n\\
    &= \bigcup_{s\in S_n} \bigcup_{\delta\in \Delta_{n}^{(s)}} \delta (\hat{s}F_n \hat{s}^{-1})\cdot \hat{s}\widehat{T}_n\\
    &= \bigcup_{s\in S_n} \bigcup_{\delta\in \Delta_{n}^{(s)}} \delta \hat{s}\cdot F_n\widehat{T}_n\\
    &= \bigcup_{\tilde{\jmath}\in \widetilde{J}_n} \tilde{\jmath}\widetilde{F}_n.
\end{align*}
It is not hard to check that, $\tilde{\jmath}\widetilde{F}_n\cap \tilde{\jmath}'\widetilde{F}_n=\emptyset$ whenever $\tilde{\jmath}\neq \tilde{\jmath}'$. Then, $(\widetilde{F}_n)_{n\in \N}$ satisfies condition (ii).
This completes the proof of the theorem.
\qed\\

The next theorem states that under stepwise monotileability, if $(F_m)_{m\in\N}$, $(T_n)_{n\in\N}$ and $(m_n)_{n\in\N}$ from the previous construction are appropriately chosen, the sequence $(F_{m_n}\widehat{T}_n)_{n\in\N}$ is congruent, exhaustive and F{\o}lner.

\begin{theorem}\label{thm:folner}
   Let $K,G,H$ be countable groups satisfying a short exact sequence as in Equation \eqref{eq:shortexact}.
    Suppose that $K$ and $H$ are both CMA groups.
    Let $(F_m)_{m\in \N}$ be a congruent and exhaustive F\o lner sequence of $K$ and $(\widehat{T}_n)_{n\in \N}$ a lifting of a  congruent and exhaustive F\o lner sequence $(T_n)_{n\in \N}$ of $H$. Consider also a  sequence $(\widehat{S}_n)_{n\in \N}$ of local centers associated to $(\widehat{T}_n)_{n\in \N}$ and an increasing sequence $(m_n)_{n\in \N}$ such that $(\widetilde{F}_n:=F_{m_n}\widehat{T}_n)_{n\in \N}$ is an exhaustive F\o lner sequence of $G$.
    If $(F_{m_n})_{n\in \N}$ is $(\widehat{S}_n)$-stepwise congruent, then $(\widetilde{F}_n)_{n\in \N}$ is also congruent.
    In particular, $G$ is a CMA group.
\end{theorem}
\begin{proof}
    Note first that the existence of an increasing sequence $(m_n)_{n\in\N}$ such that $(F_{m_n}\widehat{T}_n)_{n\in \N}$ is an exhaustive F\o lner sequence of $G$ comes from Lemma \ref{lem:subsequencefolner}, item (1).     Individually, $(F_{m_n})_{n\in \N}$ is an increasing exhaustive sequence, because $(F_m)_{m\in \N}$ is an increasing exhaustive sequence, and $(T_{n})_{n\in \N}$ is a congruent exhaustive sequence.
Therefore, since $(F_{m_n})_{n\in \N}$ is $(\widehat{S}_n)_{n\in \N}$-stepwise congruent, $(\widetilde{F}_n)_{n\in \N}$ is a congruent sequence in $G$ by Theorem \ref{thm:main}.  
The sequence $(\widetilde{F}_n)_{n\in \N}$ being congruent, exhaustive and F\o lner, we conclude that $G$ is a CMA group.
\end{proof}

\section{Examples}\label{sec:examples}

In this section we apply Theorems \ref{thm:main} and \ref{thm:folner} to several important families of group extensions, which leads us to provide new examples of CM/CMA groups. 

\subsection{Group Extensions of Locally Finite Groups} \label{subsec:locallyfinitebycma}

By Proposition \ref{prop:locallyfinite}, every locally finite group is CMA, since every exhaustive and increasing sequence of finite subgroups is congruent and also F\o lner. In what follows, we prove that all groups obtained as extensions of countable locally finite groups by CM/CMA groups are also CM/CMA. Note that if $1\longrightarrow K\longrightarrow  G\longrightarrow H\longrightarrow 1$ is a short exact sequence and $K$ is a  finite group, then $K$ is automatically $\Inn(G)$-monotileable, so the CM/CMA property for $G$ follows from the one on $H$ already by   \cite[Theorem 1.15]{salizzoni_et_al_23}.  For the more general situation where $K$ is only locally finite, we need  Theorem \ref{thm:locallyfinitebycma} that we now prove. \\

\noindent\textit{Proof of Theorem \ref{thm:locallyfinitebycma}.}
    Let G be a group that satisfies a short exact sequence as in Equation \ref{eq:shortexact}, where $K$ is locally finite and $H$ is a CMA group. Consider an increasing sequence $(A_n)_{n\in \N}$ of finite subgroups of $K$ whose union is the whole group, as in Proposition \ref{prop:locallyfinite}, and a congruent exhaustive F\o lner sequence $(T_n)_{n\in \N}$ of $H$ with associated sequence $(S_n)_{n\in\N}$ of local centers. Consider also a lifting $(\widehat{T}_n)_{n\in \N}$ of $(T_n)_{n\in \N}$ and a lifting $(\widehat{S}_n)_{n\in \N}$ of $(S_{n})_{n\in \N}$ in such a way that it contains the identity and $\widehat{T}_{n+1}=\bigsqcup_{\hat{s}\in \widehat{S}_n} \hat{s}\widehat{T}_n$ for every $n\in \N$. 
    By Lemma \ref{lem:subsequencefolner}, we know that there exists an increasing subsequence $(m_n)_{n\in \N}$ such that  $(A_{m_n}\widehat{T}_n)_{n\in \N}$ is an exhaustive F\o lner sequence of monotiles of G such that for all $n\in \N$, $\hat{s}A_{m_n}\hat{s}^{-1}\subseteq A_{m_{n+1}}$ for each $\hat{s}\in \widehat{S}_n$. In particular, $\hat{s}A_{m_n}\hat{s}^{-1}$ is a subgroup of $A_{m_{n+1}}$. Therefore, for all $\hat{s}\in \widehat{S}_n$ there exist a finite subset $\Delta_n^{(s)}\subseteq A_{m_{n+1}}$ such that $A_{m_{n+1}}=\bigsqcup_{\delta\in \Delta_{n}^{(s)}} \delta(\hat{s}A_{m_n}\hat{s}^{-1})$. This shows that the sequence $(A_{m_n})_{n\in \N}$ is $(\widehat{S}_n)$-stepwise congruent. Applying Theorem \ref{thm:folner} we conclude that $(A_{m_n}\widehat{T}_n)_{n\in \N}$ is a congruent and exhaustive F\o lner sequence in $G$. Thus, $G$ is a CMA group.\\

    Finally, if $H$ is assumed to be just CM, not necessarily amenable, and $(T_n)_{n\in\N}$ is a congruent exhaustive sequence of finite subsets of $H$, it is clear that the same construction made above will result in a congruent exhaustive sequence of finite subsets of $G$. Indeed, since the liftings $\widehat{S}_n$ are finite subsets of $G$ and $(A_n)_{n\in\N}$ is exhaustive, we can always choose $(m_n)_{n\in\N}$ such that $\hat{s}A_{m_n}\hat{s}^{-1}\subseteq A_{m_{n+1}}$ thanks to Lemma \ref{lem:subsequencefolner}, and thus $\hat{s}A_{m_n}\hat{s}^{-1}$ is a subgroup of $A_{m_{n+1}}$. The rest of the proof follows exactly as before. \\
    We conclude that if $K$ is locally finite and $H$ is CM, $G$ is CM.
\qed

\begin{example}[Group extensions of finitely supported symmetric groups]\label{ex:finitarysym}
    
Let $X$ be a nonempty set. We denote by Sym$_0(X)$ the group of all finite supported permutations of $X$, that is, the set of all bijections $\tau\colon X\to X$ such that $\{x\in X\mid \tau(x)\neq x\}$ is finite, where the group operation is given by composition. 
The group Sym$_0(X)$ is infinite, whenever $X$ is infinite, and locally finite. In the following, we consider the case where $X$ is countable.

Consider a countable group $G$ that acts by permutations on $X$. Consider the induced action $G\times \text{Sym}_0(X)\longrightarrow \text{Sym}_0(X)$ given by 
$$(g .\tau)(x):=g.\tau(g^{-1}. x ),$$

where the lower dot denotes the action of $G$ both in $X$ and $\Sym_0(X)$.

It is not difficult to verify that this defines an action of $G$ on Sym$_0(X)$ by automorphisms. Consider the semi-direct product Sym$_0(X)\rtimes G$, where $\rtimes$ is given by the previous action. A specific and well known case occurs when $X$ has a group structure. In particular, when $X=G$ and $G$ acts on itself by left multiplication. The following is immediate from Theorem \ref{thm:locallyfinitebycma}. 

\begin{corollary}\label{cor:sym}
Let $G$ be a CMA group and $X$ a countable set provided by an action $G\curvearrowright X$ by permutations. Then Sym$_0(X)\rtimes G$ is also a CMA group.
\end{corollary}

\begin{rem}
    Let $X$ be an infinite countable set. The group Sym$_0(X)$ of all finitely supported permutations of $X$ is a well-known example of a non-residually finite group. Indeed, this group contains an isomorphic copy of the simple $A_{\infty}$, the alternating group over an infinite countable set. By the same reason, Sym$_0(X)$ is not virtually solvable. 
    Therefore, Sym$_0(X)\rtimes G$ is neither residually finite nor solvable, regardless of the group $G$. Finally note that the group $\Sym_0(\Z)\rtimes \Z$ is finitely generated, thus it provides an example of a finitely generated CMA group which is neither residually finite nor virtually solvable.
\end{rem}

We now prove that for any CMA group $G$ and any inifinite countable set $X$, the group $\Sym_0(X)$ is not necessarily $\Inn(\Sym_0(X)\rtimes G)$-monotileable. More precisely, we prove that whenever the action of $G\curvearrowright X$ is transitive, $\Sym_0(X)$ is not $\Inn(\Sym_0(X)\rtimes G)$-monotileable. Let $A\subseteq X$ be a finite subset, let
$$\text{Sym}(A):=\{\tau\in \text{Sym}_0(X)\ |\ \supp(\tau)\subseteq A\},$$ 
where $\supp(\tau)=\{ x\in X\ |\ \tau(x)\neq x\}$. Note that Sym$(A)$ is a finite subgroup of Sym$_0(X)$.

 By definition, $g.\supp(\tau)=\supp(g. \tau)$ for all $g\in G$ and $\tau\in \text{Sym}_0(X)$. Therefore, one deduces that $g. \text{Sym}(A)= \text{Sym}(g. A)$ for all $g\in G$ and all finite subsets $A\subseteq X$. This identity will be useful in the proof below.

\begin{proposition}\label{prop:sym_notInnmonotileable}
    Let $G$ be a countable group and let $X$ be an infinite countable set. If $G$ acts on $X$ transitively, then $\text{Sym}_0(X)$ is not $\Inn(\Sym_0(X)\rtimes G)$-monotileable. 
\end{proposition}
\begin{proof}

We denote the group operations on $G$ and $\Sym_0(X)\rtimes G$ by central dots.
   Observe that from the definition of the group operation in $\Sym_0(X)\rtimes G$,  
    $$(\Id_X, g)\cdot (\tau, h)\cdot (\Id_X, g)^{-1}=(g.\tau, g h g^{-1})$$
    for all $\tau \in \Sym_0(X)$ and $g, h, \in G$. 
    If we identify $\Sym_0(X)$ with the subgroup $\Sym_0(X)\times \{1_G\}$, the conjugation by $(\Id_X,1_G)$ restricted to $\Sym_0(X)$ acts like 
    $$(\Id_X, g)\cdot (\tau, 1_G)\cdot  (\Id_X, g)^{-1}=(g.\tau, 1_H).$$
   This means that for all $g\in G$, the conjugation by $(\Id_X,g)$ restricted to $\Sym_0(X)$, which we denote $\varphi_g$, verifies $\varphi_g(\tau)=g.\tau$. In other words, Inn$(\Sym_0(X)\rtimes G)$ contains all $G$-shifts on $\Sym_0(X)$.\\
   Let $(F_n)_{n\in \N}$ be an increasing sequence of finite subsets of $\Sym_0(X)$ and fix $n\in \N$ such that $|F_n|\geq 2$.
   Consider two different elements $\sigma, \sigma'\in F_n$ and $z\in X$ such that $\sigma(z)\neq \sigma'(z)$. 
   Since $F_{n+1}$ is finite and each element of $F_{n+1}$ has finite support, there exists a finite subset $Y\subset X$ such that $F_{n+1}\subseteq \Sym(Y)$. 
    On the other hand, since $G$ acts on $X$ transitively, there exists $g\in G$ such that $g.z\in X\setminus Y$. 
    
Suppose that $(F_n)_{n\in \N}$ is Inn$(\Sym_0(X)\rtimes G)$-monotileable. Then, $\forall \phi\in \Inn(\Sym_0(X)\rtimes G)$, $F_{n+1}$ is a disjoint union of translations of $\phi(F_n)$. In particular, $\varphi_g(F_n)$ tiles $F_{n+1}$. As explained before, $\varphi_g(F_n)=g.F_n$, which means that there exists $J_n\subseteq \Sym_0(X)$ such that $F_{n+1}=\bigsqcup_{\tau\in J_n}\tau\circ (g.F_n)$. This implies that $\forall \tau\in J_n$, $\tau\circ (g. F_n)\subseteq F_{n+1}$. Since $g.z\not\in Y$ and $\sigma, \sigma'\in F_n$,
$$[\tau \circ (g.\sigma)](g.z)=g.z=[\tau \circ (g.\sigma')](g.z).$$
On the other hand, $(g.\sigma)(g.z)= g.\sigma(z).$ Therefore
\begin{align*}
    \sigma(z)\neq\sigma'(z) &\implies g.\sigma(z) \neq g.\sigma'(z)\\
    &\implies \tau(g.\sigma(z)) \neq \tau(g.\sigma'(z))\\
    &\iff [\tau \circ (g.\sigma)](g. z)\neq [\tau \circ (g.\sigma')](g.z),
\end{align*}
which is a contradiction.
Thus, $\tau\circ (g\cdot F_{n})\not\subseteq F_{n+1}$.
This shows that $F_{n+1}$ cannot be a disjoint union of translations of $g. F_{n}$. 
In particular, $(F_{n})_{n\in \N}$ is not $\Inn(\Sym_0(X)\rtimes G)$-monotileable. Since $(F_n)_{n\in\N}$ was arbitrarily chosen, we conclude that $\Sym_0(X)$ is not $\Inn(\Sym_0(X)\rtimes G)$-monotileable. 
\end{proof}

\end{example}

In the next example we recall the definition of Houghton's groups, and observe that they are CMA.

\begin{example}[Houghton's groups $\mathcal{H}_n$]\label{ex:houghton}
    Let $n\geq 1$ and $R_k:=\{me^{i\theta}\ |\ m\in \N, \theta=2\pi(k-1)/n \}\subset \C$
for every $1\leq k\leq n$. The disjoint union $Y_n:=\bigcup_{k=1}^n R_k$ represents $n$ discrete rays emanating from the origin in the plane. Given its combinatorial nature, we may identify $Y_n$ with the space $\{1,\dots, n\}\times \N$ and each point in $Y_n$ with a pair $(k,m)$, where $m$ is the module of the element. The Houghton's group $\mathcal{H}_n$, introduced originally in \cite{Houghton}, is defined as the group of all eventually translations of $Y_n$, {\em i.e.}, all bijections $g:Y_n\rightarrow Y_n$ for which there exists an $n$-tuple $(t_1,\dots t_n)\in \Z^n$ and a finite set $K_g\subset Y_n$ such that $g(k,m)=(k, m+t_k)$ for all $(k, m)\in Y_n\setminus K_g$. 

From its definition, the groups $\mathcal{H}_1$ and $\mathcal{H}_2$ are isomorphic to the groups $\Sym_0(\N)$ and $\Sym_0(\Z)\rtimes \Z$, respectively. Therefore, both are CMA groups, the first one because it is locally finite and the second one thanks to Corollary \ref{cor:sym}. In order to conclude the same for $n\geq 3$, we need the following two results proved by Lee in \cite{Lee}.

On the one hand, Lee showed that the derived subgroup $[\mathcal{H}_n,\mathcal{H}_n]$ of $\mathcal{H}_n$ coincides with the group $\Sym_0(Y_n)$ consisting of all permutations of $Y_n$ with finite support \cite[Lemma 2.3]{Lee}. On the other hand, assigning to each permutation $g\in \mathcal{H}_n$ its associated $n$-tuple $(t_1, \dots, t_n)$, it induces an isomorphism from the abelianization $\mathcal{H}_{n}/[\mathcal{H}_{n},\mathcal{H}_{n}]$ to the free abelian group $\Z^{n-1}$ \cite[Corollary 2.4]{Lee}. This implies the following result.

\begin{proposition}\label{prop:exacthoughton}
    For all $n\geq1$, the group $\mathcal{H}_n$ satisfies the following short exact sequence
    $$1\longrightarrow \Sym_0(Y_n)\longhookrightarrow  \mathcal{H}_n\longrightarrow \Z^{n-1}\longrightarrow 1.$$
    In particular, each $\mathcal{H}_n$ is a locally finite-by-abelian group.
\end{proposition}

From the previous proposition and Theorem \ref{thm:locallyfinitebycma} we conclude the following.

\begin{corollary}
    For all $n\geq 1$, the Houghton's group $\mathcal{H}_n$ is CMA.
\end{corollary}

\end{example}

\subsection{Wreath Product}\label{subsec:wreath}
The second family of examples we treat here consists of groups constructed by the permutational wreath product of countable groups, which generalize the usual wreath product construction.

Consider two groups $G$,$H$ and a nonempty set $X$. Suppose that $H$ acts on $X$ by permutations. 
The permutational wreath product of $G$ with $H$, with respect the action $H\curvearrowright X$, is defined as the group
$$G\wr_X H:=\left(\bigoplus_{x\in X}G \right)\rtimes H,$$
where we identify the direct sum $\bigoplus_{x\in X}G$ with the finite support functions in $G^{X}$, that is, 
$$\bigoplus_{x\in X}G=\{\varphi\in G^X: |\supp(\varphi)|<\infty\},$$
where $\supp(\varphi)=\{x\in X\ |\ \varphi(x)\neq 1_G\}$. The action $\rtimes$ is defined by $(h. \varphi)(x)=\varphi(h^{-1}. x)$.
To avoid an overcharged notation, we denote by $G^{(X)}$ the direct sum $\bigoplus_{x\in X}G$.
The identity $1_{G\wr_X H}$ of $G\wr_X H$ is the pair $(1_{G^X},1_H)$.
We remark that the classical wreath product corresponds to the case $X=H$ and $H$ acts on itself by left multiplication. In this case, we denote $G\wr H$ instead of $G\wr_H H$.

In order to remain in the class of countable groups, in what follows we shall only consider the case where $X$ is countable. 

\begin{rem}\label{rem:Grunberg}
    A wreath product $G\wr H$ is not necessarily residually finite even if $G$ and $H$ are. In fact, Grünberg proved that a wreath product $G\wr H$ is residually finite if and only if $G$ and $H$ are residually finite and $G$ is abelian or $H$ is finite \cite{Gruenberg}.
    There exists also a more general version of this result for permutational wreath products \cite{Yves}. 
    This shows that there exist wreath products that are neither residually finite nor solvable. Consider for instance the groups $\F_n\wr\F_n$, where $\F_n$ denotes the free group generated by $n$ elements, with $n\geq 2$. For an amenable example, we may consider the wreath product $A_{\infty}\wr \Z$, where $A_{\infty}$ denotes the alternating group of infinite degree.
\end{rem}

The aim of this section is to prove Theorem \ref{thm:wr}, namely that the CM/CMA property is stable under permutational wreath product. 
In order to simplify the computations below, we introduce the following notation. For $A\subseteq G^{(X)}$ and $B\subseteq H$, we denote $(A,B)=\{(a,b)\in G^{(X)}\rtimes H:a\in A, b\in B\}$.
For subsets $F\subseteq G$ and $Y\subseteq X$, we define
$$F^Y:=\{\varphi\in G^{(X)}\mid \supp(\varphi)\subseteq Y \text{ and } \forall x\in \supp(\varphi),\varphi(x)\in F \}\subseteq G^{(X)}.$$
For these sets, we identify two properties that will be useful in the sequel. First, note that $h.\supp(\varphi) =\supp(h.\varphi)$. From this it follows that $h.F^Y=F^{h. Y}$ for all $h\in H$, $F\subseteq G$ and $ Y\subseteq X$. Second, in the particular case where $F$ and $Y$ are finite sets, the number of different functions whose support is a given subset of $Y$ of size $k$ is exactly $(|F|-1)^k$, for each $0\leq k\leq |Y|$. Since there are $\binom{|Y|}{k}$ subsets of $Y$ of size $k$, this implies that $ |F^Y|=|F|^{|Y|}.$\\

We begin with a general lemma. For the statement recall from Section \ref{subsec:CMA} that  $(F_n)_{n\in \N}$ is a F\o lner sequence of a group $G$ if and only if for every finite subset $C\subseteq G$ and every $\varepsilon>0$ there exists $N\in \N$ such that for all $n\geq N$ the set $F_{n}$ is right $(C,\varepsilon)$-invariant, that is
\[ |\{g\in F_n\ |\ gC\subseteq F_n\}| \geq (1-\varepsilon)\cdot |F_n| . \]
In other words, the proportion of points of $F_n$ whose $C$-neighborhood still belongs to $F_n$ is bounded below by $(1-\varepsilon) |F_n|$. Note that if $C\subseteq C'$ are finite subsets of $G$, then 
$  \{ g\in F_n\ |\ gC'\subseteq F_n\} \subseteq  \{g\in F_n\ |\ gC\subseteq F_n\}. $

\begin{lemma} \label{lem:newfolner}
Let $(Y_n)_{n\in\N}$ be an increasing and exhaustive sequence of finite subsets of $X$ and  let ($F_n)_{n\in\N}$ be a F\o lner sequence of $G$.  Then $(F_n^{Y_n})$ is a F{\o}lner sequence of $G^{(X)}$.

\end{lemma}

\begin{proof}
    For $C\subseteq G^{(X)}$ a finite set and $\varepsilon >0$,  consider $N\in \N$ such that the support of every $c\in C$ is contained in $Y_N$. 

Since $C$ is finite, the set $C(x)=\{c(x)\mid c\in C\}$ is a finite subset of $G$ for every $x\in X$. Since $Y_N$ is finite the set $\mathcal C=\bigcup_{x\in Y_N}C(x)$ is a finite subset of $G$ as well. Now consider $M\geq N$ such that for all $m\geq M$, $F_m$ is $(\mathcal C,\widetilde \varepsilon)$-invariant for $\widetilde \varepsilon$ satisfying $(1-\widetilde\varepsilon)^{|Y_N|} \geq (1-\varepsilon)$.  

We need to estimate the proportion of functions $a\in F_M^{Y_M}$, such that $a\cdot  C$ still belongs to $ F_M^{Y_M}$. Since the elements of $C$ are supported in $Y_N$, we have that  $a\cdot  C\subset F_M^{Y_M}$ if and only if $\tilde a \cdot C\subset F_M^{Y_N}$, where $\tilde a$ denotes the restriction of $a$ to $Y_N$. By the choice of $M$, for each $x\in Y_N$ there are at least 
$(1-\widetilde\varepsilon) \cdot |F_M|$ elements in $F_M$ whose $\mathcal C$-neighborhood belongs to $F_M$. In particular, there are $\left( (1-\widetilde\varepsilon) \cdot |F_M|\right)^{|Y_N|}$ functions $\tilde a\in F_M^{Y_N}$ such that $\tilde a \cdot  C\subset F_M^{Y_N}$. Since each of these functions can be extended arbitrarily to $Y_M\setminus Y_N$, we get that there are
$$ \left( (1-\widetilde\varepsilon) \cdot |F_M|\right)^{|Y_N|} \cdot |F_M|^{|Y_M\setminus Y_N|} = (1-\widetilde\varepsilon)^{|Y_N|} |F_M|^{|Y_M|} \geq (1-\varepsilon)|F_M^{Y_M}|$$
functions whose $ C$-neighborhood is contained in $F_M^{Y_M}$. So $F_M^{Y_M}$ is $(C,\varepsilon)$-invariant as desired. The same argument shows that $F_m^{Y_m}$ is $(C,\varepsilon)$-invariant for all $m\geq M$.
\end{proof} 

\noindent\textit{Proof of Theorem \ref{thm:wr}.} Let $G\wr_XH$ be a permutational wreath product. We now show that $G\wr_XH$ is a CM/CMA group whenever $G$ and $H$ are CM/CMA groups. We start by proving the statement for CMA groups. The argument is divided depending on whether the $H$-action on $X$ has finitely or infinitely many orbits. 
\\

\noindent\textbf{Case 1: { $H\curvearrowright X$} has finitely many orbits.}
Let $G$, $H$ be CMA groups. Assume first that  $ X$ is made of a single $H$-orbit, and  let $x$ be an arbitrary point in $X$. 
Let $(F_{n})_{n\in \N}$ and $(T_n)_{n\in \N}$ be congruent and exhaustive F\o lner sequences in $G$ and $H$, respectively. Let $(J_n)_{n\in \N}$ and $(S_n)_{n\in \N}$ be sequences of local centers associated to $(F_n)_{n\in \N}$ and $(T_n)_{n\in \N},$ respectively.

In order to construct a congruent and exhaustive F\o lner sequence in $G\wr_X H$, we first build a congruent and exhaustive F\o lner sequence in $G^{(X)}$.

\begin{lemma}\label{lem:congruent}
    Let $(F_{n})_{n\in \N}$ and $(T_n)_{n\in \N}$ be congruent and exhaustive sequences in $G$ and $H$ respectively, let $x\in X$. Then, the sequence $(F_n^{T_n.x})_{n\in\N}$ is congruent, and exhaustive in $G^{(X)}$.
\end{lemma}

\begin{proof}
Let $(J_n)_{n\in \N}$ be a sequences of local centers associated to $(F_n)_{n\in \N}$. The sequence $( F_{n}^{T_n.x})_{n\in \N}$ is increasing and exhaustive because $(F_{n})_{n\in\N},$ $(T_n)_{n\in \N}$ are exhaustive increasing sequences and $X=H. x$.
Moreover, $(F_{n}^{T_n. x_i})_{n\in \N}$ is congruent because, for every $n\in\N$, the set 
$$\Gamma_n= F_{n+1}^{(T_{n+1}. x) \setminus (T_n. x)}\cdot J_n^{T_n.x}$$
satisfies the identity $F_{n+1}^{T_{n+1}.x}=\bigsqcup_{\gamma\in \Gamma_n} \gamma\cdot F_{n}^{T_{n}.x}$, which means that $(\Gamma_n)_{n\in \N}$ is a sequence of local centers for $(F_{n}^{T_n.x})$. Indeed, on the one hand, 

\begin{align*}
    \Gamma_n\cdot F_n^{T_n.x}&= (F_{n+1}^{(T_{n+1}. x) \setminus (T_n. x)}\cdot J_n^{T_n.x})\cdot F_{n}^{T_n.x}\\
    &= F_{n+1}^{(T_{n+1}. x) \setminus (T_n. x)}\cdot (J_nF_n)^{T_n.x}\\
    &=F_{n+1}^{(T_{n+1}. x) \setminus (T_n. x)}\cdot (F_{n+1})^{T_n.x}\\
    &= F_{n+1}^{T_{n+1}.x}.
\end{align*}

Hence, $F_{n+1}^{T_{n+1}.x}=\bigcup_{\gamma\in \Gamma_n} \gamma\cdot F_{n}^{T_{n}.x}$. 
On the other hand, consider two different elements $\gamma$, $\gamma'$ in $\Gamma_n$. Since they are different, they must differ at some point $z\in T_{n+1}. x$. If $z\in (T_{n+1}.x)\setminus (T_n.x)$, then 
$$\gamma(z) \cdot \varphi(z) = \gamma(z) \neq \gamma'(z) = \gamma'(z)\cdot \varphi'(z) $$
for all $\varphi,\varphi'\in F_n^{T_n.x}$. Consequently, $\gamma\cdot (F_{n}^{T_n.x})\cap \gamma'\cdot (F_{n}^{T_n.x})=\emptyset$. If now $z\in T_{n}.x$, then $\gamma(z)$ and $\gamma'(z)$ are different elements in $J_n$. So,
$\gamma(z)F_n\cap \gamma'(z)F_n=\emptyset$ and hence $\gamma\cdot (F_{n}^{T_n.x})\cap \gamma'\cdot (F_{n}^{T_n.x})=\emptyset$.\end{proof}

We now show that $G\wr_X H$ is CMA. For this recall that $(S_n)$ is the set of local centers associated to $(T_n)$. We  consider the canonical liftings $\widehat{T}_n=(1_{G^{(X)}}, T_n)$ and $\widehat{S}_n=(1_{G^{(X)}}, S_n)$. We can therefore identify $T_n$ and $S_n$ with their liftings in $G\wr_X H$, and with this identification the conjugation by elements in $\widehat{S}_n$ corresponds to the shift action of elements in $S_n$.

By Lemmas \ref{lem:newfolner} and \ref{lem:congruent}, we know that $(F_n^{T_n. x})_{n\in \N}$ is a congruent and exhaustive F\o lner sequence in $G^{(X)}$.
Thus, Lemma \ref{lem:subsequencefolner} implies  that there exists an increasing subsequence $(m_n)_{n\in \N}$ such that $[\widetilde{F}_n:=( F_{m_n}^{T_{m_n}. x}, T_n)]_{n\in \N}$ 
is F\o lner in $G\wr_X H$  and, moreover
$$F_{m_n}^{(s\cdot  T_{m_n}). x}= s.\left (  F_{m_n}^{T_{m_n}. x}\right )\subseteq F_{m_{n+1}}^{T_{m_{n+1}}. x}$$
for all $s\in S_n$, for all $n\in \N$.  We have the following lemma.

\begin{lemma}\label{lem:corona stepwise} The sequence $(F_{m_n}^{T_{m_n}. x} )_{n\in \N}$ is an $(S_n)$-stepwise congruent sequence. 
    \end{lemma} 
\begin{proof} First note that the sequence $(F_{m_n})_{n\in\N}$ is congruent and exhaustive, because $(F_n)_{n\in \N}$ is congruent and exhaustive.
Consider a sequence $(I_n)_{n\in \N}$ of local centers associated to $(F_{m_n})_{n\in\N}$.
Since $F_{m_n}^{(s\cdot T_{m_n}).x}$ is contained in $ F_{m_{n+1}}^{T_{m_{n+1}}.x}$
for all $s\in S_n$ and $n\in \N$, the operation of the group allows us conclude that the set
$$\Delta_n^{(s)}:=F_{m_{n+1}}^{({T_{m_{n+1}}. x)\setminus ((s \cdot T_{m_n}). x)}}\cdot I_n^{(s \cdot T_{m_n}). x} $$
satisfies that 
\begin{align*}
       \Delta_n^{(s)}\cdot ( F_{m_n}^{(s\cdot  T_{m_n}). x} )=&(F_{m_{n+1}}^{({T_{m_{n+1}}. x)\setminus ((s \cdot T_{m_n}). x)}}\cdot I_n^{(s \cdot T_{m_n}). x} )\cdot ( F_{m_n}^{(s\cdot  T_{m_n}). x} )\\
      =&\ F_{m_{n+1}}^{({T_{m_{n+1}}. x)\setminus ((s\cdot T_{m_n}). x)}}\cdot (I_nF_{m_n})^{(s \cdot T_{m_n}). x} \\
      =&\ F_{m_{n+1}}^{({T_{m_{n+1}}. x)\setminus ((s \cdot T_{m_n}). x)}}\cdot F_{m_{n+1}}^{(s\cdot T_{m_n}). x}\\
      =&\ F_{m_{n+1}}^{{T_{m_{n+1}}. x}}.
\end{align*} 
Moreover, arguing as in the end of proof of Lemma \ref{lem:congruent}, it follows that 
\[\delta\cdot  (F_{m_n}^{(s\cdot  T_{m_n}). x} )\cap \delta' \cdot (F_{m_n}^{(s \cdot T_{m_n}). x})=\emptyset\]
for all  $\delta\neq \delta'$ in $\Delta_n^{(s)}$.  Therefore, each $\Delta_{n}^{(s)}$ is finite, contains the identity, and satisfies $F_{m_{n+1}}^{{T_{m_{n+1}}. x}}=\bigsqcup_{\delta \in \Delta_n^{(s)}} \delta \cdot (F_{m_n}^{(s\cdot T_{m_n}). x} )$.
Thus $(F_{m_n}^{T_{m_n}. x})_{n\in \N}$ is an $(S_n)$-stepwise congruent sequence in $G^{(X)}$. \end{proof}

Therefore, by Theorem \ref{thm:folner} we conclude that $(\widetilde{F}_n)_{n\in \N}$ is a congruent exhaustive F\o lner sequence of $G\wr_X H$, hence the group $G\wr_X H$ is CMA.

\begin{rem}
   For the general case where $H\curvearrowright X$ has finitely many orbit, the argument is completely analogous to the single-orbit case. 
   If $\{x_0,\dots, x_{d-1}\}$ is a set of orbit representatives, then, as in Lemma \ref{lem:congruent}, the sequence $(\prod_{i=0}^{d-1} F_n^{T_{n}.x_i})_{n\in \N}$ is congruent and exhaustive in $G^{(X)}$ with sequence of local centers
   $$\Gamma_n=\prod_{i=0}^{d-1} (F_{n+1}^{(T_{n+1}. x_i) \setminus (T_n. x_i)}\cdot J_n^{T_n. x_i}).$$ 
   By Lemma \ref{lem:newfolner}, $(\prod_{i=0}^{d-1} F_n^{T_{n}.x_i})_{n\in \N}$ is also F\o lner.
   
   Finally, for a subsequence $(m_n)$ obtained from Lemma \ref{lem:subsequencefolner} so that $[\widetilde{F}_n:=(\prod_{i=0}^{d-1} F_{m_n}^{T_{m_n}.x_i}, T_n)]_{n\in \N}$ is F\o lner    
   and $s.(\prod_{i=0}^{d-1} F_{m_n}^{T_{m_n}.x_i})\subseteq \prod_{i=0}^{d-1} F_{m_{n+1}}^{T_{m_{n+1}}.x_i}$, we check that $(\prod_{i=0}^{d-1} F_{m_n}^{T_{m_n}.x_i})_{n\in \N}$ is $(S_n)$-stepwise congruent by considering the subsets $\Delta_n^{(s)}=\prod_{i=0}^{d-1} ( F_{m_{n+1}}^{({T_{m_{n+1}}. x_i)\setminus ((s\cdot T_{m_n}). x_i)}}\cdot I_n^{(s \cdot  T_{m_n}). x_i} )$ as in Lemma \ref{lem:corona stepwise}.
\end{rem}

\noindent\textbf{Case 2: $H\curvearrowright X$ has infinitely many orbits.} 
The arguments below  are analogous to the previous case. Let $\{x_0,x_1,\ldots\}$ be a set of orbit representatives for the $H$ action on $X$, and assume that $G$ and $H$ are CMA groups.

Let $(F_{n})_{n\in \N}$ and $(T_n)_{n\in \N}$ be congruent and exhaustive F\o lner sequences in $G$ and $H$, respectively. Let $(J_n)_{n\in \N}$ and $(S_n)_{n\in \N}$ be sequences of local centers associated to $(F_n)_{n\in \N}$ and $(T_n)_{n\in \N},$ respectively. In order to build a congruent exhaustive F\o lner sequence in $G\wr_X H$, we first consider the sequence $(U_n)_{n\in \N}$ in $G^{(X)}$ given by 
$$U_n:= \prod_{i=0}^n F_{n}^{T_n. x_i}=F_n^{\cup_{i=0}^n T_n.x_i}$$

By Lemma \ref{lem:newfolner} $(U_n)_{n\in \N}$ is a F\o lner sequence in $G^{(X)}$. Since each subset $U_n$ is defined by a finite product of subsets with pairwise disjoint supports, 
we can extend Lemma \ref{lem:congruent} to this case and conclude that $(U_n)_{n\in \N}$ is congruent and exhaustive.
Indeed, in order to prove that $(U_n)_{n\in \N}$ is congruent, it suffices to consider the sequence of local centers 
$$\Gamma_n=
    F_{n+1}^{T_{n+1}. x_{n+1}}\cdot \left (\prod_{i=0}^n (F_{n+1}^{(T_{n+1}. x_i) \setminus (T_n. x_i)}\cdot J_n^{T_n. x_i}) \right).$$
Then, remembering that $J_nF_n=F_{n+1}$, we get that  

\begin{align*}
    \Gamma_n \cdot U_n &= \left [F_{n+1}^{T_{n+1}. x_{n+1}}\cdot \left (\prod_{i=0}^n (F_{n+1}^{(T_{n+1}. x_i) \setminus (T_n. x_i)}\cdot J_n^{T_n. x_i}) \right) \right ]\cdot \prod_{i=0}^n F_n^{T_{n}.x_i} \\
    &= F_{n+1}^{T_{n+1}. x_{n+1}}\cdot \left ( \prod_{i=0}^n (F_{n+1}^{(T_{n+1}. x_i) \setminus (T_n. x_i)}\cdot (J_nF_n)^{T_n. x_i}) \right )\\
    &= F_{n+1}^{T_{n+1}. x_{n+1}}\cdot \left ( \prod_{i=0}^n (F_{n+1}^{(T_{n+1}. x_i) \setminus (T_n. x_i)}\cdot F_{n+1}^{T_n. x_i}) \right )\\
    &= F_{n+1}^{T_{n+1}. x_{n+1}}\cdot \prod_{i=0}^n F_{n+1}^{T_{n+1}.x_i}  = U_{n+1}.
\end{align*}

Moreover, as in Lemma \ref{lem:congruent}, the sets $\{\gamma\cdot U_n\}_{\gamma \in \Gamma_n}$ are pairwise disjoint, and so the sequence $U_n$ congruent (and F\o lner). 

By Lemma \ref{lem:subsequencefolner} there is an increasing subsequence $(m_n)_{n\in \N}$ such that $[\widetilde{F}_n:=( U_{m_n}, T_n)]_{n\in \N}$ 
is F\o lner in $G\wr_X H$  and, moreover
$s.( U_{m_{n}})\subseteq U_{m_{n+1}}$
for all $s\in S_n$, for all $n\in \N$.
Finally, as in the single orbit case, letting $I_n=J_{m_n+1}\cdots J_{m_{n+1}}$, we  conclude that
$$\Delta_n^{(s)}:= \prod_{i=m_n +1}^{m_{n+1}} F_{m_{n+1}}^{{T_{m_{n+1}}. x_i}}\ \cdot \prod_{i=1}^{m_{n}} ( F_{m_{n+1}}^{({T_{m_{n+1}}. x_i)\setminus ((s\cdot T_{m_n}). x_i)}}\cdot I_n^{(s \cdot  T_{m_n}). x_i} ) $$
satisfies the identity $U_{m_{n+1}}=\bigsqcup_{\delta\in \Delta_n^{(s)}} \delta\cdot (s.U_{m_n})$.
Hence,  $(U_{m_n})_{n\in \N}$ is a $(S_n)_{n\in \N}$-stepwise congruent sequence of $G^{(X)}$.
Then, Theorem \ref{thm:folner} implies that  $(\widetilde{F}_n)_{n\in \N}$ is congruent, exhaustive, and F\o lner in $G\wr_X H$. Thus, $G\wr_X H$ is a CMA group, regardless of the number of orbits.\\

Now suppose that $G$ and $H$ are CM groups, not necessarily amenable. We prove that $G\wr_X H$ is CM. Let $(F_n)_{n\in\N}$ and $(T_{n})_{n\in\N}$ be congruent exhaustive sequences of monotiles of $G$ and $H$ respectively, and
let $(J_n)_{n\in \N}$ and $(S_n)_{n\in \N}$ be sequences of local centers of $(F_{n})_{n\in \N}$ and $(T_{n})_{n\in \N}$.
Consider the sequence $(U_n)_{n\in \N}$ of subsets defined by

$$U_n:=\left \{ 
\begin{array}{ll}
    \prod_{i=0}^{d-1} F_n^{T_n.x_i} & \text{if there are finitely many orbits}. \\
    \prod_{i=0}^{n} F_n^{T_n.x_i} & \text{if there are infinitely many orbits}.
\end{array}
\right .$$
Repeating the arguments from the CMA case, we get that $U_n$ is a congruent and exhaustive sequence in $G^{(X)}$. Then, the sequences $(\widetilde{F}_n)_{n\in \N}$ defined by
$$\widetilde{F}_n:=\left \{
\begin{array}{cl}
   (\ \prod_{i=0}^{d-1} F_n^{T_n. x_i}, T_n\ ),  &  \text{if there are finitely many orbits}.\\
   (\ \prod_{i=0}^n F_n^{T_n. x_i}, T_n\ )  & \text{if there are infinitely many orbits}.
\end{array}
\right .$$
is congruent and exhaustive in $G\wr_X H$ (without having to pass to a subsequence of $(F_n)_{n\in \N}$). 
Moreover, since $\bigsqcup_{s\in S_n}sT_n=T_{n+1}$, arguing as in the CMA case  we can check that for each $s\in S_n$ and each $n\in \N$, the set
$$\Delta_n^{(s)}:=\left \{
\begin{array}{ll}
    \prod_{i=0}^{d-1} ( F_{n+1}^{(T_{{n+1}}. x_i)\setminus ((s\cdot T_{n}). x_i)}\cdot J_n^{(s\cdot T_{n}). x_i} ), & \text{if there are finitely many orbits.} \\
    F_{{n+1}}^{T_{{n+1}}. x_{i}} \cdot \prod_{i=1}^{n} ( F_{n+1}^{(T_{{n+1}}. x_i)\setminus ((s\cdot T_{n}). x_i)}\cdot J_n^{(s\cdot T_{n}). x_i} ), & \text{if there are infinitely many orbits.}
\end{array}
\right .$$
satisfies the identity $U_{n+1}=\bigsqcup_{\delta\in \Delta_n^{(s)}}\delta\cdot (s.U_{n})$.
This shows that $(U_n)_{n\in \N}$ is an $(S_n)$-stepwise congruent sequence of $G^{(X)}$ and, by Theorem \ref{thm:main}, the sequence $(\widetilde{F}_n)_{n\in \N}$ satisfies what we want.
Thus, $G\wr_X H$ is a CM group.\\

This concludes the proof of Theorem \ref{thm:wr}.
\qed\\

We immediately get the following corollary.

\begin{corollary} \label{cor:wr}
    Let $G,H$ be both CMA (resp. CM) groups. Then $G\wr H$ is also a CMA (resp. CM) group.
\end{corollary}

We end this section by noticing that $\Inn(G\wr_X H)$-monotileability is not satisfied when the action $H\curvearrowright X$ is transitive.

\begin{proposition}\label{prop:wreath_notInnmonotileable}
    Let $G, H$ be countable groups and $X$ an infinite countable set. If $H\curvearrowright X$ transitively, then the group $G^{(X)}$ is not $\Inn(G\wr_X H)$-monotileable.
\end{proposition}
\begin{proof}
The proof is very similar to that of Proposition \ref{prop:sym_notInnmonotileable}. By an analogous argument to the one presented at beginning of Proposition \ref{prop:sym_notInnmonotileable}, we have that $\Inn(G\wr_X H)$ contains all $H$-shifts on $G^{(X)}$.

   Let $(F_n)_{n\in \N}$ be an increasing sequence of finite subsets of $G^{(X)}$ and fix $n\in \N$ such that $|F_n|\geq 2$. Consider two different elements $\sigma, \sigma'\in F_n$ and $z\in X$ such that $\sigma(z)\neq \sigma'(z)$. 
   Since $F_{n+1}$ is finite and each element of $F_{n+1}$ has finite support, there exists a finite subset $Y\subset X$ such that $\supp(F_{n+1}):=\bigcup_{\tau\in F_{n+1}}\supp(\tau)\subseteq Y$.
    On the other hand, since $H$ acts on $X$ transitively, there exists $h\in H$ such that $h.z\in X\setminus Y.$ We obtain that $h.\sigma,\ h.\sigma'\in h.F_{n}$ and  
    $$(h.\sigma)(h.z)=\sigma(z)\neq \sigma'(z)=(h.\sigma')(h.z).$$
Suppose that $(F_n)_{n\in \N}$ is Inn$(G\wr_X H)$-monotileable. Then, $\forall \phi\in \Inn(G\wr_{X} H)$, $F_{n+1}$ is a disjoint union of translations of $\phi(F_n)$. In particular, $\varphi_h(F_n)$ tiles $F_{n+1}$. Since $\varphi_h(F_n)=h. F_n$, this means that there exists $J_n\subseteq G^{(X)}$ such that $F_{n+1}=\bigsqcup_{\tau\in J_n}\tau\cdot (h.F_n)$. This implies that $\forall \tau\in J_n$, $\tau\cdot (h. F_n)\subseteq F_{n+1}$. But for any $\tau\in G^{(X)}$ we have

\begin{align*}
        &[\tau\cdot(h.\sigma)](h. z)=\tau(h.z)\cdot  [h.\sigma](h. z)=\tau(h.z) \cdot\sigma(z), \text{ and}\\
        &[\tau\cdot (h.\sigma')](h. z)=\tau(h.z) \cdot [h\cdot\sigma'](h.z)=\tau(h.z) \cdot\sigma'(z).
    \end{align*}
Since $\sigma(z)\neq \sigma(z')$, this implies that $[\tau\cdot(h.\sigma)](h. z)\neq [\tau\cdot (h.\sigma')](h. z)$. Therefore, either $\sigma(z)\neq \sigma(z')\neq 1_G$ or $[\tau\cdot (h.\sigma')](h. z)\neq 1_G$, which implies that either $\supp(\tau\cdot (h. \sigma))\not\subseteq Y$ or $\supp(\tau\cdot (h.\sigma'))\not\subseteq Y$. Thus, $\tau\cdot (h.F_n)\not\subseteq F_{n+1}$, which is a contradiction. We conclude that $F_{n+1}$ cannot be a disjoint union of translated copies of $F_n$. In particular, $(F_n)_{n\in\N}$ is not $\Inn(G^{(X)})$-monotileable. Since $(F_n)_{n\in\N}$ was arbitrarily chosen, we conclude that $G^{(X)}$ is not $\Inn(G\wr_X H)$-monotileable.
\end{proof}

\subsection{Other examples}

From the results of Section \ref{subsec:locallyfinitebycma} and Section \ref{subsec:wreath} it is clear that we can obtain new examples of CMA groups by iteration of wreath products of CMA groups or group extensions of locally finite-by-CMA groups. For instance, each iterated wreath product $((\Z \wr \Z)\dots)\wr \Z)$ is CMA as well as all groups of the form $(\Sym_0(X)\rtimes G)\wr H$, provided $G$ and $H$ are CMA.
There exist also some groups defined in a similar way that do not belong to the two classes of groups mentioned above. In the following example, we treat the particular case of the {\em Lampdesigner group} (the name is borrowed from \cite{lampdesigner}).

\begin{example}[Lampdesigner groups]\label{ex:lampdesigner} 
Let $G$ and $H$ be countable groups. Consider the permutational wreath product $G\wr_H \Sym_0(H)$, defined by the canonical action of $\Sym_0(H)$ on $H$. Consider also the canonical action of $H$ on $\bigoplus_H G$ and the action of $H$ on $\Sym_0(H)$ induced by the left multiplication of $H$ on itself. 
The lampdesigner group is defined as 
$$G\boxplus H:= (G\wr_H \Sym_0(H))\rtimes H,$$
where $H$ acts on $G\wr_H \Sym_0(H))$ by 
$h\star (\tau, \sigma):= (h.\tau,h.\sigma)$.
Using the notation introduced in the previous examples, it is straightforward to check that $h\star (F^{Y}, \Sym(Z))=(F^{h\cdot  Y}, \Sym(h\cdot Z))$ for all $h\in H$, $Y,Z\subseteq H$.

Recall that $\Sym_0(H)$ is always amenable because it is locally finite. Therefore, if $G$ and $H$ are both amenable groups, then the lampdesigner group $G\boxplus H$ is also amenable. 
\begin{theorem}
    Let $G,H$ be both CM/CMA groups. Then the lampdesigner group $G\boxplus H$ is also CM/CMA.
\end{theorem}
\begin{proof}
We only consider the case where $G$ and $H$ are CMA groups.
As for the wreath product above, the proof for the CM case is analogous and is left to the reader.

Let $(F_n)_{n\in \N}$ and $(T_n)_{n\in \N}$ be congruent and exhaustive F\o lner sequences of $G$ and $H$, respectively.
Since $(T_n)_{n\in \N}$ is exhaustive, the sequence $(A_n:= \Sym(T_n))_{n\in \N}$ is a congruent and exhaustive F\o lner sequence of $\Sym_0(H)$ (see Proposition \ref{prop:locallyfinite}). Moreover, from Lemmas \ref{lem:newfolner} and \ref{lem:congruent}, we conclude that $(F_{n}^{T_n})_{n\in \N}$ is also a congruent and exhaustive F\o lner sequence of $G^{(H)}$. Therefore, from Lemma \ref{lem:subsequencefolner} there exists an increasing subsequence $(m_n)_{n\in \N}$ such that $(F_{m_n}^{T_{m_n}}, A_{n})_{n\in \N}$ is a congruent exhaustive F\o lner sequence of $G\wr_H \Sym_0(H)$.

Consider a sequence of local centers $(S_n)_{n\in \N}$ associated to $(T_n)_{n\in \N}$. 
Since $(F_{m_n}^{T_{m_n}}, A_n)_{n\in\N}$ is a congruent exhasutive F\o lner sequence, by Lemma \ref{lem:subsequencefolner} there exist an increasing subsequence $({t_n})_{n\in \N}$ such that $((F_{m_{t_n}}^{T_{m_{t_n}}}, A_{t_n}), T_n)$ is an exhaustive F\o lner sequence of monotiles of $G\boxplus H$ such that $(F_{m_{t_n}}^{s\cdot T_{m_{t_n}}},s. A_{t_n})=s\star (F_{m_{t_n}}^{ T_{m_{t_n}}}, A_{t_n})\subseteq (F_{m_{t_{n+1}}}^{ T_{m_{t_{n+1}}}},A_{t_{n+1}})$ for all $s\in S_n$, for all $n\in \N$. 
In particular, we have $s\cdot T_{m_{t_n}}\subseteq T_{m_{t_{n+1}}}$ and $s. A_{t_n}\subseteq A_{t_{n+1}}$.

Since $A_n$ is a finite group, for every $n$ and every $s\in S_n$ there exists a finite subset $\Delta_{n}^{(s)}\subset A_{t_{n+1}}$ which contains the identity and $A_{t_{n+1}}=\bigsqcup_{\delta\in \Delta_n^{(s)}} \delta\cdot (s. A_{t_n})$.
Since $\Delta_{n}(s)\subset A_{t_{n+1}}\subseteq A_{m_{t_{n+1}}}$, we also have that $\delta.(s\cdot  T_{m_{t_{n}}})\subseteq T_{m_{t_{n+1}}}$ for all $\delta\in \Delta_n^{(s)}$, $s\in S_n$, $n\in \N$.

Now, we show that $(F_{m_{t_n}}^{T_{m_{t_n}}}, A_{t_n})$ is $(S_n)$-stepwise congruent. To prove this, we must verify that $s\star(F_{m_{t_n}}^{T_{m_{t_n}}}, A_{t_n})=( F_{m_{t_n}}^{s\cdot T_{m_{t_n}}},s. A_{t_n})$ tiles $(F_{m_{t_{n+1}}}^{T_{m_{t_{n+1}}}}, A_{t_{n+1}})$ for every $s\in S_n$. 
Consider a sequence of local centers $(J_n)_{n\in \N}$ associated to $(F_{m_{t_n}})_{n\in \N}$, so that $J_nF_{m_{t_n}}=F_{m_{t_{n+1}}}$.
Then, the algebraic structure of $G\boxplus H$ allows us conclude that
\begin{align*}
    &\bigcup_{\delta\in \Delta_n^{(s)}} \left(F_{m_{t_{n+1}}}^{T_{{m_{t_{n+1}}}}\setminus \delta.(s\cdot T_{m_{t_n}})}J_n^{\delta.(s\cdot T_{m_{t_n}})},\delta \right ) \cdot \left ( F_{m_{t_n}}^{s\cdot T_{m_{t_n}}},s. A_{t_n} \right )\\
    =&\bigcup_{\delta\in \Delta_n^{(s)}} \left(F_{m_{t_{n+1}}}^{T_{{m_{t_{n+1}}}}\setminus \delta.(s \cdot T_{m_{t_{n}}})}\cdot J_n^{\delta.(s\cdot  T_{m_{t_{n}}})}\cdot F_{m_{t_n}}^{\delta.(s\cdot  T_{m_{t_{n}}})},\delta\cdot (s. A_{t_n}) \right )\\
    =&\bigcup_{\delta\in \Delta_n^{(s)}} \left(F_{m_{t_{n+1}}}^{T_{m_{t_{n+1}}}\setminus \delta.(s\cdot T_{m_{t_{n}}})} \cdot (J_nF_{m_{t_n}})^{\delta.(s\cdot T_{m_{t_{n}}})},\delta\cdot (s. A_{t_n}) \right )\\
        =&\bigcup_{\delta\in \Delta_n^{(s)}} \left(F_{m_{t_{n+1}}}^{T_{m_{t_{n+1}}}},\delta\cdot (s. A_{t_n}) \right ) =\left(F_{m_{t_{n+1}}}^{T_{m_{t_{n+1}}}},A_{t_{n+1}} \right ).
\end{align*}
Once again, a straightforward argument shows that this union is disjoint.
Thus, $(F_{m_{t_n}}^{T_{m_{t_n}}}, A_{t_n})_{n\in \N}$ is an $(S_n)$-stepwise congruent sequence. By Theorem \ref{thm:folner}, this implies that the sequence 
$$[\widetilde{F}_n:=((F_{m_{t_n}}^{T_{m_{t_n}}}, A_{t_n}), T_n)]_{n\in \N}$$
is congruent, exhaustive, and F\o lner. Therefore, $G\boxplus H$ is a CMA group. 
\end{proof}

\end{example}

\bibliographystyle{abbrv}
\bibliography{refs}

\end{document}